\documentclass[11pt]{amsart}
\usepackage{amssymb,amsmath}
\usepackage{verbatim}

\newtheorem{lemma}{Lemma}
\newtheorem*{thm}{Theorem}

\theoremstyle{remark}
\newtheorem{remark}[lemma]{Remark}

\newcommand{\col}{\,{:}\,}

\newcommand{\cB}{\mathcal{B}}

\newcommand{\cA}{\mathcal{A}}
\newcommand{\cF}{\mathcal{F}}

\newcommand{\cC}{\mathcal{C}}
\newcommand{\FF}{\mathbb{F}}
\newcommand{\F}{\FF_{\ell}}

\begin{document}



\title{An equality between two towers over cubic fields}

\author{Michael E. Zieve}
\address{
Michael E. Zieve\\
Department of Mathematics\\
Rutgers University\\
Piscataway, NJ 08854\\
USA
}

\email{zieve@math.rutgers.edu}
\urladdr{www.math.rutgers.edu/$\sim$zieve}


\begin{abstract}
Recently Bassa, Garcia and Stichtenoth constructed a tower of function
fields over $\FF_{q^3}$ having many rational places relative to their genera.
We show that, by removing the bottom field from this tower, we obtain the same
tower we would obtain by removing certain fields from a tower
constructed previously by Bezerra, Garcia and Stichtenoth.
\end{abstract}


\maketitle


Let $\ell$ be a power of the prime $p$, and let $\F$ be the field with $\ell$ elements.
If $F$ is an algebraic function field of one variable with full constant
field $\F$, we write $g(F)$ and $N(F)$ for the genus and number of degree-one
places of $F$.  We are interested in upper and lower bounds on the quantity
\[
A(\ell):=\limsup_F \frac{N(F)}{g(F)},
\]
where $F$ runs over all function fields over $\F$.  The bound
$A(\ell)\le 2\sqrt{\ell}$ follows from Weil's classical inequality
$N(F)\le \ell+1+2g(F)\sqrt{\ell}$.  Building on work of Ihara \cite{Ihara},
Drinfel'd and Vl\u{a}du\c{t} \cite{DV} improved this to $A(\ell)\le \sqrt{\ell}-1$.
Since Ihara also showed that $A(\ell)\ge\sqrt{\ell}-1$ when $\ell$ is a square,
it follows that $A(\ell)=\sqrt{\ell}-1$ for square $\ell$.

However, the value of $A(\ell)$ is not known for any nonsquare $\ell$.
Serre \cite{Se-Harvard} (see also \cite[Appendix]{KWZ2}) showed that
$A(\ell)\ge c\log\ell$ for some absolute constant $c>0$,
and subsequent authors have improved this bound in many cases.
In case $\ell=p^3$, the best known lower bound is
$A(p^3)\ge 2(p^2-1)/(p+2)$ which was proved by Zink \cite{Zink}.
Building on an example of van der Geer and van der Vlugt \cite{vdGvdV},
Bezerra, Garcia and Stichtenoth \cite{BezGS} generalized this bound to arbitrary cubic
fields, showing that $A(q^3)\ge 2(q^2-1)/(q+2)$ for every prime power $q$.
Subsequently Bassa and Stichtenoth \cite{BS}, Ihara \cite{Ihara2}, and
Bassa, Garcia and Stichtenoth \cite{BasGS} gave simpler proofs of this result.
Our purpose here is to clarify
the relationship between \cite{BasGS} and \cite{BezGS}.

The proofs of Bassa, Bezerra, Garcia and Stichtenoth are based on towers of function
fields over $\F$ (where $\ell=q^3$).  By a \emph{tower}, we mean a sequence
$\cF:=(F_1,F_2,\dots)$ of function fields over $\F$ such that
$F_1\subseteq F_2\subseteq \dots$ and $g(F_n)\to\infty$.
It is easy to see that for any tower $\cF=(F_1,F_2,\dots)$ the limit
\[
\lambda(\cF):=\lim_{n\to\infty} \frac{N(F_n)}{g(F_n)}
\]
exists and satisfies $\lambda(\cF)\le A(\ell)$.  The papers \cite{BasGS,BS,BezGS}
present three towers $\cF$ satisfying $\lambda(\cF)\ge 2(q^2-1)/(q+2)$.
Let us call a tower satisfying this bound `good'.
Bezerra, Garcia and Stichtenoth \cite{BezGS} gave two good towers $\cA$ and
$\cB$, both of which had long and complicated proofs of goodness.
These proofs were simplified in \cite{BS} and \cite{Ihara2}, but were still rather technical.
Recently, Bassa, Garcia and Stichtenoth \cite{BasGS} presented
another good tower $\cC$ having special properties allowing for a quite
simple proof of goodness.  We will show that, if we remove
the first field from $\cC$, then we obtain the same tower we would obtain by
removing certain fields from $\cB$; hence the simple argument in \cite{BasGS}
could already have been applied to the tower $\cB$ in \cite{BezGS}
(presumably the main reason this was not done in \cite{BezGS} is that
\cite{BezGS} contains an incorrect result suggesting that $\cB$ is more
complicated than $\cC$; we will explain this in Remark~\ref{err}).

We now define the towers $\cA$, $\cB$, and $\cC$.
Let $A_1=\F(a_1)$ be the rational function field, and for $n\ge 1$ let
$A_{n+1}=A_n(a_{n+1})$ where $a_{n+1}$ satisfies
\[\frac{1-a_{n+1}}{a_{n+1}^q} = \frac{a_n^q+a_n-1}{a_n}.\]
For $n\ge 1$, let $c_n$ and $b_n$ satisfy
\[
c_n^{q-1} = 1 - \frac{1}{a_n} \qquad\text{and}\qquad b_n^{q-1}=-\frac{a_n^q+a_n-1}{a_n}.
\]
Let $C_1=\F(c_1)$ and for $n\ge 1$ let $C_{n+1}=C_n(c_{n+1})$.
Let $G_1=\F(a_1)$ and for $n\ge 1$ let $H_n=G_n(b_n)$ and $G_{n+1}=H_n(a_{n+1})$.
It follows that
\[
(c_{n+1}^q-c_{n+1})^{q-1}+1 = \frac{-c_n^{q^2-q}}{(c_n^{q-1}-1)^{q-1}} \quad\text{and}\quad
(a_{n+1}b_n)^q - (a_{n+1}b_n) = -b_n
\]
for each $n\ge 1$.  By \cite[Cor.~2.2]{BezGS} and \cite[Thm.~2.2]{BasGS},
$\cA:=(A_1,A_2,\dots)$ and $\cC:=(C_1,C_2,\dots)$ are towers of function fields
over $\F$.  In \cite[Thm.~5.1]{BezGS} a similar assertion is made about
$\cB:=(G_1,H_1,G_2,H_2,\dots,G_n,H_n,\dots)$; but we will disprove one of the
claims in that result, so we will not use the result here.  Instead we note
that the following result implies at once that $\cB$ is a tower of function
fields over $\F$:

\begin{thm}
For $n\ge 2$ we have $H_n=C_n$.
\end{thm}

\begin{proof}
For $n\ge 2$ we have
\[
c_n^{q-1} = 1-\frac{1}{a_n}=-a_n^{q-1}\left(\frac{1-a_n}{a_n^q}\right)=
-a_n^{q-1}\frac{a_{n-1}^q+a_{n-1}-1}{a_{n-1}} = a_n^{q-1}b_{n-1}^{q-1},\]
so $\F(c_n)=\F(b_{n-1},a_n)$ and thus
\[
\F(c_2,\dots,c_n)=\F(b_1,a_2,b_2,a_3,\dots,b_{n-1},a_n).\]
Also
\begin{align*}
b_n^{q-1}&=-\frac{a_n^q+a_n-1}{a_n}=\frac{-a_n^q+a_n^q\frac{1-a_n}{a_n^q}}{a_n}\\
&=\frac{-a_n^q+a_n^q\frac{a_{n-1}^q+a_{n-1}-1}{a_{n-1}}}{a_n}\\
&=a_n^{q-1}\frac{a_{n-1}^q-1}{a_{n-1}}\\
&=a_n^{q-1}(a_{n-1}-1)^{q-1}\frac{a_{n-1}-1}{a_{n-1}}\\
&=a_n^{q-1}(a_{n-1}-1)^{q-1}c_{n-1}^{q-1},
\end{align*}
so $\F(c_{n-1},a_n)=\F(b_n,a_n,a_{n-1})$.
Thus for $n\ge 2$ we have
\begin{align*}
C_n&=\F(c_1,\dots,c_n)\\&=\F(c_1)\cdot\F(c_2,\dots,c_n)\\&=
\F(c_1)\cdot\F(b_1,a_2,b_2,a_3,\dots,b_{n-1},a_n).\end{align*}
Since $a_2\in C_n$ and $\F(c_1,a_2)=\F(b_2,a_2,a_1)$, it follows that
\[
C_n=\F(b_2)\cdot\F(a_1,b_1,a_2,b_2,\dots,b_{n-1},a_n).\]
Since $C_n\supseteq\F(c_{n-1},a_n)=\F(b_n,a_n,a_{n-1})$, we have $b_n\in C_n$,
so
\[
C_n=\F(a_1,b_1,a_2,b_2,\dots,a_n,b_n)=H_n.\qedhere\]
\end{proof}

\begin{remark}
If $q=2$ then $\cC=\cA$ and $H_n=G_n$ for $n\ge 1$, so by removing duplications
from $\cB$ we obtain $\cC$.  For $q>2$, one can show that $\cB$ is not a
refinement of $\cC$ (since the extension $C_2/C_1$
is not isomorphic to an extension of fields in $\cB$).
\end{remark}

\begin{remark}
The proof of our Theorem also clarifies how $\cB$ and $\cC$ relate to $\cA$.
For $n\ge 2$ we showed that $\F(c_n)=\F(b_{n-1},a_n)$ and
$\F(c_{n-1},a_n)=\F(b_n,a_n,a_{n-1})$, so $\F(c_{n-1},a_n,a_{n+1})$ contains
$b_n$ and thus contains $c_{n+1}$.  Thus $A_n\cdot C_2=C_n$, as was observed
in \cite[Rem.~8.4]{BasGS}.  Likewise the field $\F(b_{n-1},a_n,a_{n+1})$ contains
$c_n$ and thus contains $b_{n+1}$, so $A_n\cdot H_2=H_n$ (and thus $H_n=G_n$ for $n\ge 3$).
\end{remark}

\begin{remark}\label{err}
By \cite[Thm.~5.1]{BezGS}, we have $[C_{n+1}\col C_n]=q$ for $n\ge 2$.  By our
Theorem, it follows that $[H_{n+1}\col H_n]=q$ for $n\ge 2$.  This contradicts
\cite[Thm.~5.1]{BezGS} in case $q>2$, since the latter result asserts that
$[H_{n+1}\col H_n]=q^2-q$.  The mistake in the proof of \cite[Thm.~5.1]{BezGS} becomes
clear upon inspection, since they observe that $H_1/G_1$ is visibly a Kummer
extension of degree $q-1$, and later assert without proof that a similar
argument implies $[H_n\col G_n]=q-1$ for all $n\ge 2$ (whereas we showed above that
$H_n=G_n$ for $n\ge 3$).
\end{remark}


\end{document}